\documentclass{amsart}
\usepackage{amssymb,amscd,amsthm,verbatim}

\newtheorem{proposition}{Proposition}
\newtheorem{theorem}[proposition]{Theorem}
\newtheorem{corollary}[proposition]{Corollary}
\newtheorem{lemma}[proposition]{Lemma}

\newtheorem{question}[proposition]{Question}

\newcommand{\integers}{\mathbb Z}
\newcommand{\reals}{\mathbb R}
\newcommand{\complexes}{\mathbb C}
\newcommand{\sphere}{\mathbb S}
\newcommand{\Proj}{\mathbb P}

\newcommand{\redhom}{\tilde{H}}
\newcommand{\Ind}{\mathrm{Ind}}
\newcommand{\A}{{\mathcal A}}
\newcommand{\PP}{{\mathcal P}}
\newcommand{\FF}{{\mathcal F}}
\newcommand{\del}{\partial}

\def\ignore#1{\relax}

\begin{document}

\title[The sign representation for Shephard groups]
{The sign representation for Shephard groups}

\author[Peter Orlik]{Peter Orlik}
\email{orlik@math.wisc.edu}
\address{Dept.\@  of Mathematics\\
         University of Wisconsin\\
         Madison, WI 53706}

\author[Victor Reiner]{Victor Reiner}
\email{reiner@math.umn.edu}
\address{School of Mathematics\\
         University of Minnesota\\
         Minneapolis, MN 55455}

\author[Anne V. Shepler]{Anne V. Shepler \\
\vspace{2ex}
{{\tiny Dedicated to Louis Solomon on his seventieth birthday}}
}
\email{ashepler@math.ucsc.edu}
\address{Dept.\@ of Mathematics\\
         University of California at Santa Cruz\\
         Santa Cruz, CA 96054}

\keywords{Coxeter group, unitary reflection group, Shephard group, 
regular complex polytope, arrangement of hyperplanes, Milnor fiber}

\thanks{Work of second author partially supported by 
NSF grant DMS-9877047.  Work of third author partially supported by
NSF grant DMS-9971099.\\
\indent Email contact: {\tt ashepler@math.ucsc.edu.}}

\begin{abstract}
Shephard groups are unitary reflection
groups arising as the symmetries of  regular complex
polytopes.  For a Shephard group, we identify the representation carried by
the principal ideal in the coinvariant algebra
generated by the image of the product of all linear forms
defining reflecting hyperplanes.  This representation
turns out to have many equivalent guises making it
analogous to the sign representation of a finite Coxeter group.
One of these guises is (up to a twist) the 
cohomology of the Milnor fiber for the isolated singularity at $0$
in the hypersurface defined by any homogeneous invariant of minimal degree.

\end{abstract}

\maketitle

\section{Introduction}
\label{introduction}

Let $W$ be a finite reflection group acting in a Euclidean
space $V$, that is, a finite subgroup of $GL(V)$ generated by reflections.
An important role in
the structure and representation theory of $W$ is played
by its {\it sign character} 
$$
\epsilon: W \rightarrow \integers^{\times}=\{ \pm 1 \}.
$$
This character appears in many different guises:
\begin{enumerate}
\item[$\bullet$] $\epsilon(w) = \det(w)=\det^{-1}(w).$
\item[$\bullet$] $\epsilon$ is the character of $W$ acting
on the top (reduced) cohomology group of the unit sphere $\sphere^{\dim(V)-1}$.
\item[$\bullet$] If $R$ is a set of Coxeter generators for $W$,
then $\epsilon$ is the  virtual character 
$\sum_{J \subset R} (-1)^{|R-J|}\,  \Ind_{W_J}^W 1_{W_J}$,
where $W_J$ is the parabolic subgroup generated by $J$, and
$1_{W_J}$ denotes its trivial character.
\item[$\bullet$] Let $S=\complexes[V]$ denote the ring of polynomial
functions on $V$, and $I$ the ideal generated by the $W$-invariant
polynomials of positive degree.  Then the quotient ring $S/I$ is
a graded ring which is finite-dimensional over $\complexes$, and
whose nonvanishing graded component of top degree $(S/I)_{t}$ 
carries the representation $\epsilon$.
\item[$\bullet$] This top graded component $(S/I)_{t}$ 
can also be described as the principal ideal $Q\cdot (S/I)$
within $S/I$, where $Q$ is the product of the linear forms 
defining the reflecting hyperplanes for $W$.
\end{enumerate}

This paper concerns an analogue of $\epsilon$ for
the class of unitary reflection groups known as Shephard groups.
Let $V$ be a finite-dimensional complex unitary space.
Recall that a {\it unitary reflection group} is a finite subgroup
$G \subset GL(V)$ 
generated by {\it unitary reflections},
i.e., elements of finite order that fix a hyperplane in $V$.  Such
groups include the finite Euclidean reflection groups,
called {\em Coxeter groups}, 
and were classified by Shephard and Todd \cite{ShephardTodd}.  

{\it Shephard groups} are 
the symmetry groups of the {\it regular complex polytopes} defined
and classified by Shephard \cite{Shephard} (see also Coxeter \cite{Coxeter}).
These groups generalize the finite reflection groups
which occur as the Euclidean symmetry groups
of regular convex polytopes,
or equivalently, those whose Coxeter diagrams are unbranched.
In particular, each Shephard group can be generated 
by a distinguished set $R$ of $\ell:= \dim V$ generators 
which yield a particularly nice presentation for the group, in fact,
a presentation which can be expressed by 
a ``Coxeter-like'' diagram which is unbranched; 
see Section~\ref{Shephard-groups}.

We show in this paper  that every Shephard group $G$
has a representation, defined over $\integers$, which occurs in
many guises, analogous to $\epsilon$.
We introduce some notation to make this precise.
Most of our notation follows \cite{OrlikTerao}.

Let $\PP$ be a regular complex polytope in a unitary space $V$ of dimension $\ell$
having $G$ as its group of unitary automorphisms.
Let $\Delta$ be the {\it order complex}
of its poset of proper faces, that is, the simplicial complex of totally
ordered subsets in this poset. Let $R$ be a distinguished set of generators
for $G$ (as defined in Section \ref{Shephard-groups} below).  
For $J \subset R$, let 
$G_J$ denote the subgroup generated by $J$.  

As before, let $S$ denote the algebra $\complexes[V]$ of polynomial
functions on $V$, and $I$ the ideal generated by the $G$-invariants
of positive degree.  Let $d$ denote the minimal degree of a $G$-invariant,
and let $f_1$ denote any homogeneous $G$-invariant of this degree
(this turns out to define $f_1$ uniquely up to 
a scalar multiple; see Lemma~\ref{case-by-case}).
The {\it Milnor fiber} of the singularity at $0$ on the hypersurface $f_1^{-1}(0)$ 
is the level set $F:=f_1^{-1}(1)$, where we regard 
$f_1$ as a map $f_1:V \rightarrow \complexes $.
Let $K$ denote the ideal of $S$ generated by the first partial derivatives
$\frac{\del f_1}{\del x_1},\ldots,\frac{\del f_1}{\del x_\ell}$.
Let $Q$ denote the product $\prod_{H \in \A}{\alpha_H}$ where
$\A$ is the collection of reflecting hyperplanes $H$, and 
$\alpha_H$ is any linear form that vanishes
on $H$.  Given a graded vector space $U= \oplus_{i} U_i$ carrying a
graded representation of $G$ with character $\chi_{U_i}:G \rightarrow \complexes$
on $U_i$, define its {\it graded character}
$$
\chi_{U,t}(g) : = \sum_{i} \chi_{U_i}(g) t^i.
$$

After establishing notation and reviewing facts about unitary reflection
groups in Section \ref{review} and Shephard groups in Section~\ref{Shephard-groups},
we prove the following result, which is essentially a collection of 
previously known results.

\begin{theorem}
\label{collection}
Let $G$ be a Shephard group with distinguished generators $R$.
Then the following graded (complex) representations of $G$ are equivalent:
\begin{enumerate}
\item[(i)] $S/K\otimes \det^{-1}$.
\item[(ii)] the representation $U$ affording the graded character
$$
\chi_{U,t}(g) = \frac{\det(1-gt^{d-1})}{\det(g-t)}.
$$
\end{enumerate}
Furthermore, as ungraded representations, both are equivalent to the
following representations defined over $\integers$:
\begin{enumerate}
\item[(iii)] \label{virtual-character} the dual
of the virtual representation
$\sum_{J \subset R} (-1)^{|R-J|}\,  \Ind_{G_J}^G 1_{G_J}$.
\item[(iv)] the representation on the (reduced) cohomology 
$\redhom^{\ell-1}(F,\complexes)$ of the Milnor fiber.
\item[(v)] the representation on the (reduced) cohomology 
$\redhom^{\ell-1}(\Delta, \complexes)$.
\end{enumerate}
\end{theorem}

Our main result, proven in Section~\ref{main-section}, describes another
natural occurrence of this representation.  Let $\phi: S \rightarrow S/I$ be
the composite map
$$
S \overset{Q}{\rightarrow} S \rightarrow S/I
$$
where the first map is multiplication by $Q$ and the second is the
canonical surjection.
\begin{theorem}
\label{Qmap-isomorphism}
For any Shephard group $G$, the kernel of $\phi$
is the ideal $K$ generated by the first partial derivatives of $f_1$.
Therefore, $S/K$ maps isomorphically onto the principal ideal $Q \cdot (S/I)$
within $S/I$.
\end{theorem}

\noindent
It is known (see Lemma \ref{relative-invariants} below) that
$$
g(Q) = det^{-1}(g)Q \text{ for all }g\text{ in }G.
$$
Consequently, Theorem~\ref{Qmap-isomorphism} shows that
the graded representation carried by $Q \cdot (S/I)$ is 
equivalent (up to a shift in grading) with those in (i), (ii) of
Theorem~\ref{collection}, and equivalent as an ungraded representation
with those in (iii), (iv), and (v) of Theorem~\ref{collection}. 

Note that in a suitable coordinate system, the polynomial
$f_1 := x_1^2 + \ldots + x_{\ell}^2$
is a minimal degree invariant for
any Coxeter group.
Every Coxeter group $W$ 
acting on $\reals^\ell$
can be considered as a unitary reflection group
acting on $\complexes^\ell$.  
The sphere $\sphere^{\ell-1}:= \{v \in \reals^\ell: f_1(v) = 1 \}$ 
is a $W$-equivariant strong deformation retract of
the Milnor fiber $F := \{v \in \complexes^\ell: f_1(v) = 1 \}$; see
e.g.~\cite{Orlik2}.
Thus, when $G$ is simultaneously a Coxeter group and a Shephard group,
the sign character described in the introduction and the
sign representation described in Theorems~\ref{collection}
and \ref{Qmap-isomorphism} coincide.  
  
Section~\ref{remarks} contains remarks and open questions.

\section{Notation and review of unitary reflection groups}
\label{review}

Let $V$ be an $\ell$-dimensional {\it unitary space}, that is, a 
$\complexes$-vector space of dimension $\ell$ with a positive definite Hermitian form.
A {\it unitary reflection} (or {\it pseudo-reflection})
is a non-identity element $g$ of $GL(V)$ of finite order 
which fixes some hyperplane $H$ in $V$, called the
{\it reflecting hyperplane} for $g$.  A finite
subgroup $G \subset GL(V)$ is called a {\it unitary group generated by reflections}
(or {\it u.g.g.r.}) if it is generated by unitary reflections.  A u.g.g.r.\@
is {\it irreducible} if $V$ contains no $G$-invariant subspaces.
Irreducible u.g.g.r.'s were classified by 
Shephard and Todd \cite{ShephardTodd}.
They proved
that u.g.g.r's are distinguished by a rich invariant theory
which we discuss next.  
Good references for most of this material are
\cite{Humphreys, OrlikTerao, Stanley1}.

A subgroup $G \subset GL(V)$ acts on the dual space $V^*$ in the usual
(contragradient) way: for $g$ in $G$,  $f$ in $V^*$, and $v \in V$ we have
$$
g(f)(v) = f( g^{-1}(v)).
$$
This extends to an action on the symmetric algebra $S:=Sym(V^*)$, which we 
can 
view as the algebra of polynomial functions $f:V \rightarrow \complexes$.
When $G$ is finite, it is well-known \cite[\S 1]{Stanley1} that the 
subalgebra of invariant  polynomials $S^G$ is a finitely generated 
$\complexes$-algebra, and $S$ is a finitely generated module over $S^G$.
Shephard and Todd \cite{ShephardTodd} and Chevalley \cite{Chevalley}
proved the following.

\begin{theorem}
Let $G \subset GL(V)$ be a finite subgroup. Then
$S^G$ is isomorphic to a polynomial algebra generated by $\ell$ algebraically independent
homogeneous elements $f_1,\ldots,f_\ell$ if and only if $G$ is a u.g.g.r.\@ $\qed$
\end{theorem}

We call $f_1, \ldots, f_\ell$ a set of {\it basic invariants} of $G$.
Let $I$ be the ideal in $S$ generated by $f_1, \ldots, f_\ell$.
Note that although the invariants $f_1,\ldots,f_\ell$ themselves 
are not unique,
their degrees and the ideal $I$ are uniquely determined by $G$.
Let $d$ be the minimum non-zero degree in $I$, and
note that there 
can be more than one $G$-invariant of degree $d$ up to scaling.  
On the other hand, this does not
happen when $G$ is a Shephard group; see Lemma~\ref{case-by-case}.

The fact that $S$ is finite over $S^G$ means that $f_1,\ldots,f_\ell$ form
a {\it homogeneous system of parameters} ({\it h.s.o.p.})
for $S$, and consequently also an $S$-{\it regular sequence}, 
since $S$ is {\it Cohen-Macaulay};
see e.g.\@  \cite[\S 3]{Stanley1} and \cite{Stanley2}.
This says that $S/I$ is a graded {\it complete intersection}, and rings
of this form satisfy a version of Poincar\'e duality:

\begin{lemma} 
\label{complete-intersection-lemma}
Let $S/L$ be a graded complete intersection, that is,
$L$ is an ideal in $S$ generated by an
h.s.o.p.\@  $h_1,\ldots,h_\ell$
which is also an $S$-regular sequence.
Let $t_i:=\deg(h_i)$.  Then
\begin{enumerate} 
\item[(i)]\cite[\S 8, \S 9]{Stanley1} \cite[p.~12]{Stanley2}
$S$ is a Gorenstein ring of Krull dimension $0$, with top non-zero degree
$$
\tau:=\sum_{i=1}^{\ell} (t_i-1).
$$
Consequently, the bilinear pairing
$$
(S/L)_j \times (S/L)_{\tau-j} \rightarrow (S/L)_{\tau} \cong \complexes
$$
is non-degenerate.
\item[(ii)]\cite[p.~187]{SchejaStorch} $(S/L)_{\tau}$ is spanned by the image of the
Jacobian determinant 
$$
Jac(h_1,\ldots,h_\ell):=
  \det\left(\frac{\del h_i}{\del x_j}\right)_{i,j=1,\ldots,\ell}.
\qed
$$
\end{enumerate}
\end{lemma}

For a u.g.g.r.\@ $G$,
the algebra $S/I$ is called the {\it coinvariant algebra}.
A theorem of Chevalley \cite{Chevalley} asserts that
$S/I$ is equivalent to the regular representation
as an ungraded $G$-representation.
We wish to be
explicit about the occurrences of certain degree one characters in
this representation.  Let $\A$ denote the collection of
reflecting hyperplanes
of the unitary reflections in $G$, and for each such unitary reflection,
let $\alpha_H$ be a linear form
that vanishes on its reflecting
hyperplane $H$. Let $e_H$ denote the order of the cyclic subgroup of $G$
which fixes $H$.  Given any degree one character 
$\chi: G \rightarrow \complexes^{\times}$, for each hyperplane $H \in \A$
there is a unique integer $e_{H,\chi}$ with $0 \leq e_{H,\chi} < e_H$
defined by $\chi(g) = \det(g)^{-e_{H,\chi}}$ for all unitary reflections $g$
fixing $H$.  One can then define a (minimal) $\chi$-{\it relative invariant} 
$$
Q_{\chi}:= \prod_{H \in \A} (\alpha_H)^{e_{H,\chi}},
$$
i.e. $g (Q_{\chi}) = \chi(g) Q_{\chi}$ for all $g$ in $G$.

\begin{lemma} \cite[Proposition 4.12]{Stanley1}
\label{relative-invariants}
The set of $\chi$-relative invariants $S^{G,\chi}$ is
a free $S^G$-module of rank $1$ with $Q_{\chi}$ as generator. $\qed$
\end{lemma}

The following particular cases of $Q_{\chi}$ are important for what follows:
$$
\begin{aligned}
J(G) &:=Q_{\det}=\prod_{H \in \A} (\alpha_H)^{e_H-1},\\
H(G) &:=Q_{\det^2}=\prod_{H \in \A} (\alpha_H)^{e_H-2},\\
Q(G) &:=Q_{\det^{-1}}\prod_{H \in \A} \alpha_H.
\end{aligned}
$$
When no confusion will result, we will use $J, H, Q$ to refer
to $J(G), H(G), Q(G)$, respectively.
Note that one has by definition
\begin{equation}
\label{QH=J}
QH=J.
\end{equation}
Let $t$ denote the top degree of
$S/I$ (see Lemma~\ref{complete-intersection-lemma}).

\begin{lemma}\cite{Steinberg} \cite[p.~283]{ShephardTodd}
\label{JnotinI}
For any u.g.g.r.\@ $G$, $J$ is equal to the Jacobian determinant 
$Jac(f_1,\ldots,f_\ell)$ up to a scalar multiple.
Consequently, the image of $J$ spans $(S/I)_{t}$,
and in particular, $J$ does not lie in $I$. \qed
\end{lemma}

\section{Shephard groups}
\label{Shephard-groups}

We now turn to the special case of Shephard groups, which enjoy
special properties not shared by all u.g.g.r.'s.
For a more detailed treatment of Shephard groups,  see
Coxeter's wonderful book \cite{Coxeter}.

A {\it regular complex polytope} $\PP$ in $V$ is a collection of
complex affine subspaces of $V$, called {\it faces} of $\PP$,
satisfying certain conditions \cite[p.~115]{Coxeter}. 
One of these conditions is that the group $G \subset GL(V)$
of unitary automorphisms of $\PP$ acts transitively on the
maximal flags of faces in $\PP$.  Such a group $G$ is called
a {\it Shephard group}, and will always be an irreducible u.g.g.r.
One obtains a {\em distinguished} 
set of generators $R:=\{r_0,\ldots,r_{\ell-1}\}$
for a Shephard group $G$ as follows:  let 
\begin{equation}
\label{base-flag}
{\mathcal F_0:}=( F_0 \subset F_1 \subset \cdots \subset F_{\ell-2} \subset F_{\ell-1})
\end{equation}
be a fixed maximal flag of (proper) faces in $\PP$, which we will
call the {\it base flag}.
For each $i$, choose $r_i$ to be a generator for the (cyclic) group that
stabilizes (not necessarily pointwise) each
$F_j$ with $j \neq i$.  Let $p_i$ denote the
order of $r_i$; then there exist positive integers $q_i \geq 3$ such
that $G$  has the following very simple presentation with respect to
these generators:
$$
\begin{aligned}
 r_i^{p_i} &= 1, \\ 
 r_i r_j &= r_j r_i \text{ if }|i-j|>1, \\
\underset{q_i\text{ letters}}{\underbrace{r_i r_{i+1} r_i r_{i+1} \cdots}}
 &  = 
\underset{q_i\text{ letters}}{\underbrace{ r_{i+1} r_i r_{i+1} r_i \cdots}}.
\end{aligned}
$$
The Shephard group $G$ with the above presentation is denoted by
the shorthand {\it symbol}
\begin{equation}
\label{shorthand}
p_0[q_0]p_1[q_1]p_2 \cdots p_{\ell-2}[q_{\ell-2}]p_{\ell-1}.
\end{equation}
It may also be represented by a ``Coxeter-like'' linear diagram with vertices
labeled by the $p_i$ and edges labeled by the $q_i$.
The classification of Shephard groups is relatively short.
There is one infinite family $r[4]2[3]2[3]\cdots2[3]2$ isomorphic
to the wreath product $C_r \wr S_\ell$ of a cyclic group with a
symmetric group (corresponding to $G(r,1,\ell)$ in the notation of
Shephard and Todd \cite{ShephardTodd}).  There is  a finite list
of exceptional Shephard groups:
\begin{enumerate}
\item[$\bullet$]  symmetry groups of real regular polytopes 
(Coxeter groups with unbranched diagrams)
\item[$\bullet$] $p_0[q]p_1$ 
where $p_0, p_1 \geq 2$ and $q\geq 3$ satisfy
$p_0 = p_1$ if $q$ is odd,
$$
 \frac{1}{p_0} + \frac{1}{p_1} + \frac{2}{q}  > 1,
$$
and at least one of $p_0, p_1$ is $>2$. There are twelve such groups.
\item[$\bullet$] $2[4]3[3]3$
\item[$\bullet$] $3[3]3[3]3$
\item[$\bullet$] $3[3]3[3]3[3]3$.
\end{enumerate}

From the invariant-theoretic point of view, Shephard groups
have the following extra properties, which were verified using
the above classification in 
\cite[Corollaries 5.4 and 5.8, Theorem~5.10]{OrlikSolomon4}:

\begin{lemma}
Let $G$ be a Shephard group, and $d$ the minimal degree
of a $G$-invariant.  
\begin{enumerate}
\label{case-by-case}
\item[(i)] Up to scalar multiples, there is a unique $G$-invariant $f_1$ of degree $d$.
\item[(ii)] The hypersurface $f_1^{-1}(0)$ has an isolated critical
point at $0$ in $V$, i.e., $f_1$ defines a smooth hypersurface
in the projective space $\Proj(V)$.
\item[(iii)]  The Hessian determinant
$$
\begin{aligned}
Hess(f_1)&:=Jac\left( \frac{\del f_1}{\del x_1}, 
                  \ldots, \frac{\del f_1}{\del x_\ell} \right)\\
     &= \det\left(\frac{\del^2 f_1}
                       {\del x_i x_j} \right)_{i,j=1,\ldots,r}
\end{aligned}
$$
is equal (up to scalar multiple) to
$H(G) =Q_{\det^2}= \prod_{H \in \A}(\alpha_H)^{e_H-2}. \qed$
\end{enumerate}
\end{lemma}
\noindent
We remark that property (iii) above actually {\it characterizes}
the union of the Coxeter and the Shephard groups among all u.g.g.r.'s, 
see \cite[Theorem 6.121]{OrlikTerao}.

We next discuss some consequences of Lemma \ref{case-by-case}.
As in the introduction, for a Shephard group $G$, let $K$ denote
the ideal in $S$ generated by the first partial derivatives 
$\frac{\del f_1}{\del x_1},\ldots,\frac{\del f_1}{\del x_\ell}$.

\begin{corollary}
\label{S/K-Gorenstein}
For a Shephard group $G$ with notation as above,
$S/K$ is a graded complete intersection with top
degree $\ell(d-2)$ and with $(S/K)_{\ell(d-2)}$ spanned by the image of $H$.
Consequently, for any $f$ in $S-K$, there exists $f'$ in $S$
with $\overline{ff'} = \overline{H}$ in $S/K$. 
\end{corollary}
\begin{proof}
The fact that $0$ is an isolated singular point of $f_1^{-1}(0)$
implies that the first partial derivatives of $f_1$
form an $S$-regular sequence by 
\cite[Chapter~V, p.~137, Exercise~5]{Bourbaki}.
The rest follows from Lemma~\ref{complete-intersection-lemma} and
Lemma~\ref{relative-invariants}.
\end{proof}

Part (ii) of Lemma~\ref{case-by-case} also has the following
consequences for the topology of the {\it Milnor fiber} $F:=f_1^{-1}(1)$
(see \cite{OrlikSolomon1, OrlikSolomon2} and the references
therein). 

\begin{theorem}
\label{Fary-Milnor}
Let $G$ be a Shephard group with minimal 
degree invariant $f_1$ of degree $d$ as
above.  
Then
\begin{enumerate}
\item[(i)] The Milnor fiber $F$ is homotopy equivalent to a wedge of $(d-1)^\ell$
spheres of dimension $\ell-1$.
\item[(ii)] There is a $G$-equivariant isomorphism 
$S/K \rightarrow \redhom^{\ell-1}(F;\complexes)$. $\qed$
\end{enumerate}
\end{theorem}
 
Let $\Delta$ denote the {\it order complex}
of the poset of proper faces of the regular complex polytope $\PP$.
In other words, $\Delta$ is the simplicial complex having vertex
set indexed by the proper faces of $\PP$ and simplices
corresponding to flags of nested faces.   Note that the choice
of a base flag $\FF_0$ as in definition
\eqref{base-flag} then corresponds
to the choice of maximal face in $\Delta$ which we call the
{\it base chamber}.

The following is proven in \cite[Thms. 4.1 and 5.1]{Orlik2}
somewhat nonconstructively; see \cite{OrlikSolomon6, Szydlik} for more
explicit case-by-case constructions that use the classification
of Shephard groups.

\begin{theorem}
\label{Milnor-fiber-complex}
The geometric realization of $\Delta$ is $G$-equivariantly
isomorphic to a ($G$-equivariant) strong 
deformation retraction of the Milnor fiber $F:=f_1^{-1}(1)$. $\qed$
\end{theorem}

The following corollary is \cite[Corollary 5.3]{Orlik2};
see \cite{Bjorner1, Bjorner2} for more on Cohen-Macaulay complexes.

\begin{corollary}
\label{Cohen-Macaulayness}
$\Delta$ is  a {\it Cohen-Macaulay complex}.
\end{corollary}

\begin{proof}
Theorems~\ref{Fary-Milnor} and
\ref{Milnor-fiber-complex} imply that $\Delta$ has only top-dimensional
cohomology.  The same follows for all links of faces in $\Delta$,
since these are always joins of order complexes
of regular complex polytopes which are {\it medial polytopes} of $\PP$
(see \cite[p.~116]{Coxeter} or the proof of
Lemma~\ref{stabilizers} below).

\end{proof}

We wish to describe explicitly the permutation action of $G$ on faces
of $\Delta$.  Most of the following lemma seems implicit in the discussion
of medial polytopes from \cite[p.~116]{Coxeter}, but we include 
a proof
for the sake of completeness.

\begin{lemma}
\label{stabilizers}
Let $\PP$ be a regular complex polytope, and $G$ its Shephard group.
Then $G$ acts simply transitively on maximal flags of faces
in $\PP$, and hence on maximal faces (chambers) of $\Delta$.

More generally, consider a partial flag
$$
\FF=( F_{a_0} \subset F_{a_1} \subset \cdots \subset F_{a_{k}} ) 
$$
contained in the base flag $\FF_0$ of definition \eqref{base-flag},
or equivalently, a face contained in the base chamber of $\Delta$.
Then the stabilizer subgroup 
within $G$ of $\FF$ is the subgroup $G_J$ generated by 
the subset of distinguished 
generators $J :  = R - \{r_{a_0}, r_{a_1}, \ldots, r_{a_{k}}\}$.
\end{lemma}
\begin{proof}
The first assertion is \cite[p.~116, lines 1-2]{Coxeter}.

For the second assertion, 
note that $G_J$ is a subset of the stabilizer of $\FF$,
and hence it suffices to show that they have the same cardinality.
By the first assertion, the order of a Shephard group $G$ is the 
number of maximal flags in the corresponding polytope $\PP$,
and a group element $g$ may be identified with the image $g\FF_0$
of the base flag $\FF_0$.
In particular, the stabilizer of $\FF$ has the same cardinality
as the set of maximal flags in $\PP$ which pass through the partial
flag $\FF$.  This cardinality is clearly the product of the
numbers of maximal flags in each interval
$$
[F_{a_{i-1}}, F_{a_{i}}]:= \{ \text{faces $F$ in $\PP$ with }
F_{a_{i-1}} \subset F \subset F_{a_{i}} \}
$$
for $i=0,\ldots,k+1$ 
(where we adopt the convention that 
$a_{-1} := -1, a_{k+1} := \ell$,
$F_{-1}: = \emptyset$,
and $F_{\ell}:=V$).  However, each such interval is again
the poset of faces in a regular complex polytope $\PP_i$,
the {\it medial polytope} \cite[p. ~116]{Coxeter} associated 
with $F_{a_{i-1}} \subset F_{a_{i}}$.  Since the Shephard group associated
to $\PP_i$ may be identified with the subgroup $G_{J_i}$ where 
$$
J_i: = \{r_{a_{i-1}+1}, r_{a_{i-1}+2}, \ldots, r_{a_i-2},r_{a_i-1} \},
$$
we conclude that the stabilizer of $\FF$ has cardinality 
$\prod_{i=0}^{k+1} |G_{J_i}|$.
On the other hand, since $J = \cup_{i=0}^{k+1} J_i$, and $G_{J_r}, G_{J_s}$
commute for $r \neq s$ by the presentation of $G$ discussed
in Section \ref{Shephard-groups}, we 
conclude that $|G_J| = \prod_{i=0}^{k+1} |G_{J_i}|$,
as desired.
\end{proof}

We are now in a position to prove Theorem~\ref{collection}.
\vskip .2in
\noindent
{\it Proof of Theorem~\ref{collection}.}
The equivalence of (i) and (ii) is, up to a twist by $\det^{-1}$,
exactly \cite[Theorem 1.3]{OrlikSolomon2}.
The equivalence of (i) and (iv) is the main theorem of \cite{OrlikSolomon1}.
The equivalence of (iv) and (v) follows from
 Theorems \ref{Fary-Milnor} and
\ref{Milnor-fiber-complex}.

The equivalence of (iii) and (v) will follow from 
the equivalent statement that the virtual representation
$\sum_{J \subset R} (-1)^{|R-J|}\, \Ind_{G_J}^G 1_{G_J}$
is equivalent to the reduced homology representation
$\redhom_{\ell-1}(\Delta, \complexes)$.
We use the Hopf trace formula:
\begin{equation}
\label{Hopf-trace} 
\sum_{i \geq -1} (-1)^i Trace\left( g|_{\tilde{C}_i(\Delta, \complexes)}\right) =
\sum_{i \geq -1} (-1)^i Trace\left( g|_{\redhom_i(\Delta, \complexes)}\right),
\end{equation}
where $\tilde{C}_i(\Delta,\complexes)$ denotes the $i$-dimensional chain group
in the augmented simplicial chain complex that computes the
reduced homology $\redhom_\cdot(\Delta, \complexes)$.
Lemma~\ref{stabilizers} implies that the permutation action of $G$ 
on $\tilde{C}_i(\Delta)$ is the direct sum of the coset actions
$G/G_J$ as $J$ ranges over subsets of $R$ with $|R-J|=i+1$,
so that its character is the sum of induced characters 
$\Ind_{G_J}^G 1_{G_J}$ over the same set of $J$'s.  Thus the left-hand-side
of \eqref{Hopf-trace} is $(-1)^{|R|-1}$ times the
character $\sum_{J \subset R} (-1)^{|R-J|}\, \Ind_{G_J}^G 1_{G_J}$.  
Meanwhile, Corollary~\ref{Cohen-Macaulayness} implies that all 
$\redhom_i(\Delta,\complexes)$ vanish except when $i=\ell-1$, so
that the right-hand side of \eqref{Hopf-trace} is $(-1)^{\ell-1}$ times
the character of the homology representation $\redhom_{\ell-1}(\Delta, \complexes)$.  $\qed$

\section{Proof of Theorem~\ref{Qmap-isomorphism}}
\label{main-section}

Before proving Theorem~\ref{Qmap-isomorphism},
we review some facts about anti-invariant forms.
The action of a group $G \subset GL(V)$ on $S$ induces an action on the set
of derivations (or vector fields) on $V$, $\text{Der}_S \simeq S \otimes V$. 
This in turn induces an action on the set of differential $1$-forms on $V$, 
$\Omega^1 := \text{Hom(Der}_S, S) \simeq S \otimes V^*$.
The set $\text{Der}_S$ is a free $S$-module
with basis  $\{ \frac{\del}{\del x_i} \}$,
and $\Omega^1$ is a free $S$-module with
basis $\{ dx_i \}$.
The modules $\text{Der}_S$ and $\Omega^1$ 
inherit gradings from $S$:
we say that a derivation or form has degree $p$
if the coefficient of each  $\frac{\del}{\del x_i}$
or $dx_i$ is homogeneous of degree $p$.

When $G$ is a u.g.g.r.,  
the set of invariant derivations is a free $S^G$-module
of dimension $\ell$, and we call a set of homogeneous generators 
{\em basic derivations}.  Basic derivations are not
determined uniquely by $G$, but their degrees,
called the {\em coexponents} of $G$, are.  
The coexponents are intimately connected 
with the invariant theory of $G$.
If $G$ acts irreducibly in $V$ 
and hence in $V^*$, then the representation $V^*$ occurs in $S/I$ 
with multiplicity $\ell$ and 
in homogeneous components given by the coexponents.
See \cite[Chapter 6]{OrlikTerao} 
for more on invariant theory and coexponents, especially
for Shephard groups.

A differential $1$-form $\omega$ is called 
{\em anti-invariant} if it is relatively invariant 
with respect to the $\det^{-1}$ character of $G$, 
i.e., 
$$
 g(\omega) = {\det}^{-1} (g) \ \omega ,
$$
for any $g$ in $G$.
Let $(\Omega^1)^{\det^{-1}}$ be the space of anti-invariant $1$-forms.
This space is a free $S^G$-module of rank $\ell$.

We construct generators for $(\Omega^1)^{\det^{-1}}$ from a set of basic 
derivations.  Let $n_1, \ldots, n_\ell$ be the coexponents for
$G$ and let $\theta_1, \ldots, \theta_\ell$ be a set of basic derivations
with $\deg(\theta_i) = n_i$.
We follow \cite{SheplerTerao}.
Let $\Omega^1(\A)$ be the $S$-module
of {\em logarithmic $1$-forms} with poles along $\A$ 
(see \cite{OrlikTerao}):
$$
\Omega^{1}(\A) := \left\{ \frac{\eta}{Q} :  \ \eta\in \Omega^{1},\  
d\left(\frac{\eta}{Q}\right) \in \frac{1}{Q}  \Omega^{2} \right\}
$$
where $d$ is the exterior differentiation and
$\Omega^{2} := \Omega^1 \bigwedge_S \Omega^1$
is the $S$-module of differential $2$-forms.
By \cite[Theorem 6.59]{OrlikTerao}, 
$\theta_1, \ldots, \theta_{\ell}$ is an $S$-basis
for the module $D(\A)$ of $\A$-derivations, 
$$
D(\A) := \{\theta\in {\rm Der_{S}} \mid \theta(Q) \in QS\}.
$$
By the contraction 
$\left< ~, ~\right>$ of a $1$-form and a derivation,
the $S$-modules
$D(\A)$ and $\Omega^1 (\A)$ are $S$-dual to
each other \cite[Theorem 4.75]{OrlikTerao} .   
Let
$\{\omega_{1}, \ldots , \omega_{\ell}  \}
\subset \Omega^{1}(\A) $
be the basis of $\Omega^{1}(\A)$ dual to
$\{ \theta_{1}, \ldots , \theta_{\ell}  \}$:
$\omega_i$ is the unique element of $\Omega^{1}(\A)$
satisfying  
$\left<\theta_{i}, \omega_{j}  \right> =
\delta_{ij} $  (Kronecker's delta).
The group $G$ acts naturally on 
$\Omega^{1}(\A)$, and 
each $\omega_i$ is invariant since 
the $\theta_1, \ldots, \theta_\ell$ are invariant.
Let $\mu_i := Q \omega_i$
for each $i$.  Then since $Q$ is anti-invariant,
each $\mu_i$ is an anti-invariant $1$-form
of degree $(\deg Q - n_i)$.
The following lemma is an application of \cite[Thm.~ 1]{SheplerTerao}
or \cite[Prop.~1]{Shepler}.
\begin{lemma}
\label{anti-invariantforms}
Let $G$ be a u.g.g.r. Then
$$
 Q\ \Omega^1(\A) = (\Omega^1)^{\det^{-1} } \otimes_{S^G} S,
$$
and hence
$$
 (\Omega^1)^{det^{-1}} = S^G \mu_1 \oplus \cdots \oplus S^G \mu_\ell.
$$
\end{lemma}
The differential forms $\omega_i$,
and hence the $\mu_i$, may be constructed explicitly as follows.
Let $M$ be the coefficient matrix 
of $\{ \theta_1, \ldots, \theta_\ell \}$, i.e.,
the matrix whose $(i,j)$ entry is the polynomial coefficient
of $\frac{\del}{\del x_j}$ in $\theta_i$.
Using Saito's criterion, Terao showed that 
$\det(M) = Q$; see  \cite[Chapter 6]{OrlikTerao}.
Let $w_{ij}$ be the $(i,j)$ entry of $M^{-1}$.  Then
each $w_{ij}$ is a rational function with denominator $Q$. 
For each $i$, $\omega_i$ is the rational differential form
$$
 \omega_i := \sum_{j=1}^\ell  \ w_{ji} \ dx_j.
$$

\begin{lemma}
\label{QPinI}
Let $G$ be any u.g.g.r.\@ and $f \in S$ a $G$-invariant of positive degree.
Then $Q \frac{\partial f}{\partial x_i}$ lies in $I$ for $i=1,2, \ldots,l$.
In particular, if $f_1$ is a $G$-invariant of minimal positive degree,
and $K$ is the ideal generated by its first partial derivatives, then
$$Q K \subset I$$
and hence $K \subset \ker \phi$.
\end{lemma}
\begin{proof}
Let $df$ be the exterior derivative of $f$.  Since $f$ is invariant,
$df$ is invariant, and hence 
$Q df$ is an anti-invariant $1$-form. 
By Lemma~\ref{anti-invariantforms}, $Q df$ can thus
be written as a combination of the $\mu_i$ with coefficients
from $S^G$:
\begin{equation}
\label{Qdf}
  Qdf = h_1 \mu_1 + \ldots + h_\ell \mu_\ell. 
\end{equation}
Since
$\deg(Qdf) = \deg Q + \deg f -1 > \deg Q - n_i = \deg \mu_i$ 
for each $i$, each $h_i$ must have positive degree
and thus lie in $I$.
By comparing the coefficient of $d x_i$ on each side of equation \eqref{Qdf}
above,
we see that each $Q \frac{\del f}{\del x_i}$
is in $I$.
\end{proof}

\noindent
{\it Proof of Theorem~\ref{Qmap-isomorphism}.}
By Lemma \ref{QPinI}, we only need to show that
$\ker \phi \subset K$.  Assume for the sake of contradiction
that $f$ is in $\ker \phi$,
i.e. $Qf$ lies in $I$, but $f$ is not in $K$. 
By Lemma~\ref{S/K-Gorenstein}, there exists some $f' \in S$ 
with $\overline{ff'} = \overline{H}$ in $S/K$.
Consequently we have 
$$
\begin{aligned}
H &= ff'+ k \text{ for some }k\text{ in }K \\
QH & = Qff' + Qk.\\
\end{aligned}
$$
The left-hand side is $J$ by equation \eqref{QH=J}.
The right-hand side lies in $I$, because $QK \subset I$ by Lemma~\ref{QPinI}
and $Qf \in I$ by assumption.  We conclude that $J$ lies in $I$,
contradicting Lemma~\ref{JnotinI}.
$\qed$.

\section{Remarks and questions}
\label{remarks}

We conclude with some remarks and open questions.

\subsection{Coxeter complexes for u.g.g.r.'s}
Lemma \ref{stabilizers} shows that for any Shephard group $G$
and a distinguished set of generators $R$,
the simplicial complex $\Delta$ has an alternate construction
that parallels the construction of {\it Coxeter complexes} for
Coxeter groups \cite[\S 1.15]{Humphreys}: 
it is the unique simplicial complex whose poset
of faces is isomorphic to the poset of cosets of ``standard parabolic''
subgroups 
$$
P(G,R):= \{gG_J\}_{\substack{g \in G, \\J \subset R}}
$$
ordered by reverse inclusion.

This construction may be carried out more generally.  Given any pair
$(G,R)$ of a group $G$ and a finite set of generators $R$ which is
minimal with respect to inclusion, one can form the poset of
cosets $P(G,R)$ as above.  It is not always true that this is
the face poset of an abstract simplicial complex.  However,
it {\it is} always the face poset of 
a {\it regular cell complex} $\Delta(G,R)$,
in which all faces are isomorphic to simplices, but the
intersection of a pair of faces need not be a face of 
each; see \cite{Bjorner2, Stanley3} for more on such cell complexes.

The cell complex $\Delta(G,R)$ shares many of the pleasant properties of Coxeter
complexes, and its homology carries representations of the group $G$.
We always have the Hopf trace formula relating the virtual character
from Theorem~\ref{collection} (iii) 
to the alternating sum of the
$G$-representations on the homology groups of $\Delta(G,R)$.
Unlike the Coxeter or Shephard group cases, this homology
need not be concentrated in a single dimension, so that this virtual
character need not be a genuine character (even up to sign).

\begin{question}
For u.g.g.r.'s $G$ other than Coxeter or Shephard groups,
do there exist minimal generating sets $R$ for which the
``Coxeter complex'' $\Delta(G,R)$ carries homology
representations related to the homology representations of
Milnor fibers $f^{-1}(1)$ for some interesting $G$-invariant or
relative-invariant $f$?
\end{question}

We have investigated this a tiny bit for the presentations
of u.g.g.r.'s given in \cite{BroueMalleRouquier}
with inconclusive results.

\subsection{Shellability}
Corollary~\ref{Cohen-Macaulayness} suggests the following question.

\begin{question}
Is $\Delta$ {\it shellable} for any regular complex polytope $\PP$?
More strongly, is the poset of faces of $\PP$ lexicographically shellable?
\end{question}

\noindent
See \cite{Bjorner2} for the definitions of shellable
and lexicographically shellable, and 
the relation to being Cohen-Macaulay.

Shellability of $\Delta$ is well-known for the Shephard groups which
are also Coxeter groups, so one might try to resolve this question 
by appealing to the 
classification of the remaining Shephard groups.  

For the infinite family associated
with the Shephard group $G(r,1,\ell)$, one of the two associated
regular complex polytopes is Shephard's {\it generalized cross-polytope}
$\beta^r_\ell$ \cite{Shephard, Coxeter}.  Its face poset (with top
element removed) is easily seen to be the 
$\ell$-fold Cartesian product of a poset having $r+1$ elements
with one bottom element and the rest atoms.  This makes it very
easy to produce a lexicographic shelling, answering both questions affirmatively.

For Shephard groups with $\ell=2$,
$\Delta$ is a connected graph, hence trivially shellable.
We have not checked whether the face poset of
$\PP$ is lexicographically shellable.

For the remaining three exceptional cases, we have checked
using the computer algebra package GAP that the method
used by Solomon and Tits to shell Tits buildings 
\cite[Chapter~IV,
\S~6]{Brown}
(ordering the maximal faces by any linear ordering that respects
distance from the base chamber) seems not to give a shelling.  
We also have no candidate for a lexicographic shelling of the
face poset of $\PP$.

\subsection{Retractions}
Orlik and Solomon \cite{OrlikSolomon4, OrlikSolomon5}
observed interesting and mysterious connections 
between the invariant theory
for a Shephard group $G$ having symbol
$p_0[q_0]p_1[q_1]p_2 \cdots p_{\ell-2}[q_{\ell-2}]p_{\ell-1}$
and the ``associated'' Coxeter group $W$ having symbol
$2[q_0]2[q_1]2 \cdots 2[q_{\ell-2}]2$.  We hypothesize further
connections between the Coxeter complexes of $G$ and $W$.

Let $\Delta_G$ be the simplicial complex $\Delta$ which was
associated to $G$ in Section \ref{review}, and let $\Delta_W$
be the corresponding complex associated to $W$ (that is, the
{\it Coxeter complex} of $W$). 
In the special case of the infinite family of Shephard
groups $G=G(r,1,\ell)$, it is not hard to see that there are
many well-defined simplicial inclusions and retractions
$$
\begin{aligned}
\iota:\Delta_W &\hookrightarrow\Delta_G\\
\rho:\Delta_G &\twoheadrightarrow \Delta_W
\end{aligned}
$$
satisfying $\rho \circ \iota = \mathrm{id}_{\Delta_W}$
which preserve the ``coloring'' of vertices by the
distinguished generators of each group (this coloring assigns
the color $i$ to vertices of $\Delta$ which correspond to 
$i$-dimensional faces in the regular polytope, or equivalently,
to those which correspond to cosets of the form $gG_{R-\{r_i\}}$).

We give an example of such an inclusion-retraction
pair using Shephard's generalized cross-polytope $\beta^r_\ell$ 
and the usual cross-polytope $\beta^2_\ell$ as the 
regular complex polytopes associated to $G$ and $W$, respectively.
In this case, the maps
can be determined from their restriction to $0$-faces in
the associated polytopes.
A typical $0$-face in the generalized cross-polytope is $\omega^{k}e_i$
with $\omega = e^{\frac{2\pi i}{r}}$, $0 \leq k \leq r-1$, $1 \leq i \leq \ell$.
Define
$$
\begin{aligned}
\iota(+e_i) &= e_i = \omega^{0} e_i, \\
\iota(-e_i) &= \omega^1 e_i,\\
 & \\
\rho(\omega^{k}e_i) &= 
\begin{cases}
+e_i &\text{ if }k=0 \\
-e_i &\text{ if }k>0.
\end{cases}
\end{aligned}
$$

\begin{question}
\label{retraction-question}
Can similar inclusions and retractions be defined for any Shephard group $G$
and associated Coxeter group $W$?
\end{question}

One motivation for this question comes from the Weyl group 
$W$ of a {\it finite reductive group}
$G'$ (see \cite[Introduction]{BroueMalleMichel}).  Such groups $G'$ have a
$BN$-pair structure which gives rise to a simplicial complex $\Delta_{G',B,N}$
known as a {\it Tits building} \cite{Brown, Garrett}.  The Tits building has many
subcomplexes isomorphic to $\Delta_W$ (called {\it apartments})
and there is a {\it canonical retraction} 
$\Delta_{G,B,N} \twoheadrightarrow \Delta_W$
onto any apartment that preserves the natural coloring of vertices by
the Coxeter generators of $W$.  A positive answer to Question \ref{retraction-question}
would provide further support for the following analogy:
$$
\begin{matrix}
\text{Weyl group }W&\text{Shephard group }G &\text{finite reductive group }G'\\
\text{Coxeter system }(W,R)&\text{``Shephard system'' }(G,R)
         &BN\text{-pair }(G',B,N)\\
\text{sign character} & 
\text{Theorem~\ref{collection}'s representation} &\text{Steinberg representation}\\
\text{Coxeter complex }\Delta_W&\text{``Coxeter complex'' }\Delta_G&
\text{Tits building }\Delta_{G',B,N}
\end{matrix}
$$
Furthermore, whenever a retraction as in Question \ref{retraction-question} exists,
one can fill in the question mark in the
following diagram of simplicial retractions
$$
\begin{CD}
    ?    @>>>         \Delta_{G} \\
  @VVV                   @VVV\\
 \Delta_{G',B,N}  @>>> \Delta_{W} \\
\end{CD}
$$
with a simplicial complex defined using the usual pullback construction.
We hypothesize that this pullback complex plays the role of the
Tits building for the yet-to-be-defined {\it spetses} investigated by Brou\'e,
Malle, and Michel \cite{BroueMalleMichel}.  By analogy to groups with
$BN$-pair, perhaps one can define the spetses to be 
the group of vertex-color-preserving simplicial automorphisms of this complex?

\section{Acknowledgments}
The second author would like to thank Michel Brou\'e for an inspiring
series of talks at the University of Minnesota on unitary reflection
groups, and for helpful conversations.

\newcommand{\journalname}[1]{\textrm{#1}}
\newcommand{\booktitle}[1]{\textrm{#1}}

\end{document}